\documentclass[oneside,a4paper,11pt,reqno]{amsart}
\usepackage{amssymb,amsmath,amsthm}

\newtheorem{hypo}{Hypothesis}

\newtheorem{thm}[hypo]{Theorem}

\def\I{\mathcal{I}}

\def\PP{\mathbb{P}}
\def\RR{\mathbb{R}}
\def\ZZ{\mathbb{Z}}

\title[Self-intersections of Lorentz process]{An asymptotic estimate of the variance of the 
self-intersections of a planar periodic Lorentz process}

\date\today

\author{Fran\c{c}oise P\`ene}
\address{Universit\'e de Brest,
UMR CNRS 6205, Laboratoire de Math\'ematique de Bretagne Atlantique,
6 avenue Le Gorgeu, 29238 Brest cedex, France}
\email{francoise.pene@univ-brest.fr}

\keywords{Sinai billiard, Lorentz process, self intersection \\
Fran\c{c}oise P\`ene is supported by the french ANR projects GEODE (ANR-10-JCJC-0108) and
PERTURBATIONS (ANR-10-BLAN-0106)}
\begin{document}

\begin{abstract}
We consider a $\mathbb Z^2$-periodic planar Lorentz process with strictly convex obstacles
and finite horizon. This process describes the
displacement of a particle moving in the plane with unit speed and with elastic
reflection on the obstacles.
We call number of self-intersections of this Lorentz process the number $V_n$ of couples
of integers $(k,\ell)$ smaller than $n$ such that the particle hits a same obstacle 
both at the $k$th and at the $\ell$th collision times.
The aim of this article is to prove that the variance of $V_n$ is equivalent to $n^2$
(such a result has recently been proved for simple planar random walks in \cite{DU}).
\end{abstract}
\maketitle
\section{Introduction}
\label{}
We consider a finite number of convex open sets $O_1,...,O_I\subset \RR^2$ with
boundary $C^3$-smooth and with non null curvature. We repeat these sets $\mathbb Z^2$-periodically by defining $U_{i,\ell}=O_i+\ell$ for every $(i,\ell)\in\{1,...,I\}\times\mathbb Z^2$.
We suppose that the closures of the $U_{i,\ell}$ are pairwise disjoint.
We assume that {\bf the horizon is finite}, which means that every line meets the boundary of 
at least one obstacle (i.e. there is no infinite free flight).
We consider a particle moving in the domain $Q:=\mathbb R^2\setminus\bigcup_{i=1}^I\bigcup_{\ell\in\mathbb Z^2}U_{i,\ell}$ with unit speed and with respect to
the Descartes reflection law at its reflection times (reflected angle=incident angle).
We assume that the particle starts from $[0,1[^2\cap Q $ 
with uniform distribution in position and in speed.
The {\bf Lorentz process} describes the evolution of the particle in $Q$.
Because of the $\mathbb Z^2$-periodicity, it
is strongly related to
the Sinai billiard, the ergodic properties of which have been studied namely by Sinai
in \cite{Sinai70} (for its ergodicity), Bunimovich and Sinai \cite{BS80,BS81}, Bunimovich, Chernov and Sinai
\cite{BCS90,BCS91} (for central limit theorems), Young \cite{Young98}
(for exponential rate of decorrelation).
The similarity of behaviour of the Lorentz process with a simple planar random walk
has been investigated by many authors (\cite{SzV,DSzV},...).
The number of auto-intersections up to time $n$ of a random walk $(\tilde S_n)_n$ is
$\tilde V_n:=\sum_{k,\ell=1}^n\mathbf{1}_{\tilde S_k=\tilde S_\ell}$. 
This quantity is linked with random walks in random sceneries \cite{Bol,DU}. Recently,
in \cite{DU}, Deligiannidis and Utev proved that $Var(\tilde V_n)\sim\tilde cn^{2}$ with an explicit $\tilde c$. This improved the estimation in 
$O(n^2\log n)$ by Bolthausen \cite{Bol}.
For the Lorentz process, we define $(\I_k,S_k)$ in $\{1,...,I\}\times\mathbb Z^2$
for the index of the obstacle hit at
the $k$-th reflection time
($(I_0,S_0)$ being the index of the obstacle at the reflection time just before time 0).
Recall that $(k^{-1/2}S_k)_{k\ge 1}$ admits an asymptotic positive variance matrix $\Sigma^2$.
We call {\bf number of self-intersections} of the Lorentz process up to the $n$-th reflection time the quantity
$V_n:=\sum_{k,\ell=1}^n{\mathbf 1}_{S_k=S_\ell,\I_k=\I_\ell} $.
In \cite{FP09b}, we proved that $\mathbb E[V_n]\sim c_0 n\log n$ with
$c_0:=\frac{\sum_{i=1}^I(|\partial O_i|^2)}{(\sum_{i=1}^I|\partial O_i|)^2\pi\sqrt{\det \Sigma^2}}, $
where $|\partial O_i|$ stands for the length of $\partial O_i$.
In \cite{FP09b}, $Var(V_n)=O(n^2\log n)$ was enough for our study of the planar Lorentz process in random scenery.
Our proof of the following result uses decorrelation and 
precised local limit theorems established in \cite{FP09b}. 
It provides an alternative strategy to the one of \cite{DU}.
\begin{thm}\label{THM}
$Var(V_n)\sim c n^{2}$ 
with
$$c:=c_0^2\left(1+2J-\frac{\pi^2}6\right)\mbox{ and } 
J:=\int_{[0,1]^3}\frac{(1-(u+v+w)){\bf 1}_{\{u+v+w\le 1\}}\,  du\, dv\, dw}
    {uv+uw+vw}.$$
\end{thm}
\section{Proof of Theorem \ref{THM}}
Observe that the distribution of $(S_k-S_0,\I_k)_k$ under $\PP$ and under $\bar\nu$ considered
in \cite{FP09b} are the same (by $\mathbb Z^2$-periodicity and by
construction of $\bar\nu$).
We write $E_{k,\ell}:=\{S_k=S_\ell,\I_k=\I_\ell\}$. According to \cite{FP09b}, we have
\begin{equation}\label{TLL}
\PP(E_{k,\ell})=\PP(E_{0,|\ell-k}|)=
   c_1|\ell-k|^{-1}+O(|\ell-k|^{-2}),\ \mbox{with}\ \   
c_1:=\frac{\sum_{i}^I(\PP(\I_0=i))^2}{2\pi\sqrt{\det \Sigma^2}}=\frac{c_0}2.
\end{equation}
Observe that we have
$V_n=n+2\sum_{1\le k<\ell \le n}{\mathbf 1}_{S_k=S_\ell,\I_k=\I_\ell}$ 
and so
$$
Var(V_n)=4\sum_{1\le k_1<\ell_1 \le n}\sum_{1\le k_2<\ell_2 \le n}D_{k_1,\ell_1,k_2,\ell_2}
= 8A_1+8A_2+8A_3+4A_4,
$$
with $D_{k_1,\ell_1,k_2,\ell_2}:=\PP(E_{k_1,\ell_1}\cap E_{k_2,\ell_2})
       -\PP(E_{k_1,\ell_1})\PP( E_{k_2,\ell_2})$ and 
$$A_1:=\sum_{1\le k_1<\ell_1\le k_2<\ell_2\le n}D_{k_1,\ell_1,k_2,\ell_2},\ 
A_2:=\sum_{1\le k_1\le k_2<\ell_1\le\ell_2\le n}D_{k_1,\ell_1,k_2,\ell_2},$$ 
$$A_3:=\sum_{1\le k_1< k_2<\ell_2<\ell_1\le n}D_{k_1,\ell_1,k_2,\ell_2}, \ 
A_4:=\sum_{1\le k<\ell\le n}[\PP(E_{k,\ell})
       -(\PP(E_{k,\ell}))^2].$$
\begin{itemize}
\item \underline{Control of $A_1$}.

Due to \cite{FP09b}, if $k_1<\ell_1\le k_2<\ell_2$, then
$|D_{k_1,\ell_1,k_2,\ell_2}|\le C_1 \tau_1^{k_2-\ell_1}/((\ell_1-k_1)
         (\ell_2-k_2))  $
for some $C_1>0$ and some $\tau_1\in(0,1)$. Hence $A_1= O(n\log^2 n)=o(n^2)$.
\item \underline{Control of $A_4$}.

Due to (\ref{TLL}) or \cite{SzV}, $A_4\le C_2\sum_{1\le k<\ell\le n}(\ell-k)^{-1}=O(n\log n)=o(n^2)$.
\item \underline{Control of $A_2$}.

According to \cite{FP09b}, we have
\begin{equation}\label{A2}
A_2=\!\!\!\!\!\!\!\!\!\!\!\!\!\!\!\!\sum_{1\le k_1< k_2<\ell_1<\ell_2\le n}\!\!\!\!\!\!\! c_1^2
    \left[\left(\sum_x \frac{e^{-\frac{\langle (\Sigma^2)^{-1}x,x\rangle}{2}\left(
         \frac 1{k_2-k_1}+\frac 1{\ell_1-k_2}+\frac 1{\ell_2-\ell_1}\right)}}
   {2\pi\sqrt{\det \Sigma^2}(k_2-k_1)(\ell_1-k_2)(\ell_2-\ell_1)}\right)
  -\frac{1}{(\ell_1-k_1)(\ell_2-k_2)}\right]+o(n^2),
\end{equation}
where $\sum_x=\sum_{x\in\mathbb Z^2\, :\, |x|\le ||S_1||_\infty\min(k_2-k_1,\ell_1-k_2,\ell_2-\ell_1)}$ and $\langle\cdot,\cdot\rangle$ is the usual scalar product in $\RR^2$.
\begin{itemize}
\item First $\displaystyle
 A_{2,0}:=\sum_{1\le k_1< k_2<\ell_1<\ell_2\le n}
    \sum_x \frac{e^{-\frac{\langle (\Sigma^2)^{-1}x,x\rangle}{2}\left(
         \frac 1{k_2-k_1}+\frac 1{\ell_1-k_2}+\frac 1{\ell_2-\ell_1}\right)}}
   {2\pi\sqrt{\det \Sigma^2}(k_2-k_1)(\ell_1-k_2)(\ell_2-\ell_1)}$
$$=\sum_{(k_1,m_0,m_1,m_2)\in E_n}
    \sum_{|x|\le||S_1||_\infty \min(m_0,m_1,m_2)}
   \frac{e^{-\frac{\langle (\Sigma^2)^{-1}x,x\rangle}{2}\left(
         \frac 1{m_0}+\frac 1{m_1}+\frac 1{m_2}\right)}}
   {2\pi\sqrt{\det \Sigma^2}m_0m_1m_2},$$
with $E_n:=\{(k_1,m_0,m_1,m_2)\in \ZZ_+\, :\, 
          k_1+m_0+m_1+m_2\le n\}$. Observe that, using a comparison series-integral, we obtain
\begin{equation}\label{clef}
\sup_{||S_1||_\infty \le a\le 3||S_1||_\infty}
  \left|\sum_{x\in\mathbb Z^2\ :\ |x|\le am}e^{-\frac{\langle(\Sigma^2)^{-1} x,x\rangle}{2m}} -
 2\pi m\sqrt{\det\Sigma^2}\right|=O(\sqrt{m}).
\end{equation}
$$\mbox{So}\ \ \ \ \ \ \ \ \ \ \ A_{2,0} =\sum_{(k_1,m_0,m_1,m_2)\in E_n}
\frac{1+O(\min(m_0,m_1,m_2)^{-1/2})}{m_0m_1+m_0m_2+m_1m_2}\sim n^2 J.$$
\item Second 
$\sum_{1\le k_1< k_2<\ell_1<\ell_2\le n}\frac{1}{(\ell_1-k_1)(\ell_2-k_2)}=A_{2,1}+2A_{2,2},$
$$\mbox{with}\ A_{2,1}:=\sum_{k=1}^n\sum_{\max(1,2k-n)\le m\le k}\frac{n-(2k-m)+1}{k^2}
\le \sum_{k=1}^n\sum_{m=0}^k\frac n{k^2}=O(n\log n)=o(n^2),$$

$\displaystyle \mbox{and}\ \ 
 A_{2,2}:=\sum_{1\le k<\ell\le n}\sum_{\max(0,k+\ell-n)\le m\le k}\frac{n-(k+\ell-m)+1}{k\ell}$
\begin{eqnarray*}
&=&\sum_{\ell=1}^{\lfloor{n/2}\rfloor}\sum_{k=1}^{\ell-1}\sum_{m=0}^k\cdots+
  \sum_{\ell=\lfloor{n/2}\rfloor+1}^n\sum_{k=1}^{n-\ell}\sum_{m=0}^k\cdots+
\sum_{\ell=\lfloor{n/2}\rfloor+1}^n\sum_{k=n-\ell+1}^{\ell-1}\sum_{m=k+\ell-n}^k\cdots\\
&=&o(n^2)+\sum_{\ell=1}^{\lfloor{n/2}\rfloor}\sum_{k=1}^{\ell-1}\frac{2(n-\ell)-k}{2\ell}+
  \sum_{\ell=\lfloor{n/2}\rfloor+1}^n\sum_{k=1}^{n-\ell}\frac{2(n-\ell)-k}{2\ell}+
\sum_{\ell=\lfloor{n/2}\rfloor+1}^n\sum_{k=n-\ell+1}^{\ell-1}\frac{(n-\ell)^2}{2k\ell}\\
&=&o(n^2)+\sum_{\ell=1}^{\lfloor{n/2}\rfloor}\frac{4n-5\ell}{4}+
  \sum_{\ell=\lfloor{n/2}\rfloor+1}^n\frac{3(n-\ell)^2}{4\ell}+
\sum_{\ell=\lfloor{n/2}\rfloor+1}^n\frac{(n-\ell)^2}{2\ell}\log\left(\frac \ell{n-\ell}\right)\\
&\sim& n^2\left(-\frac 1 8 +\frac 34\log 2+\frac I2\right),
\end{eqnarray*}
with
\begin{eqnarray*}
I&:=&\int_{1/2}^1
     \frac{(1-u)^2} {u}\log\left(\frac u{1-u}\right)\, du\\
&=&\left[ Li_2(u)+\frac 12\left(u+\log u(u^2+\log u-4 u)
+\log(1-u)(-u^2+4 u-3)
\right)\right]_{1/2}^1,
\end{eqnarray*}
with $Li_2(z):=\sum_{k\ge 1}\frac{z^k}{k^2}$. So $I=Li_2(1)-Li_2(1/2)+\frac 14-\frac{\log^2 2}2
-\frac 32 \log 2
=\frac{\pi^2}6-(\frac{\pi^2}{12}-\frac{\log^2 2}2)+\frac 14-\frac{\log^2 2}2
-\frac 32 \log 2=\frac{\pi^2}{12}+\frac 14-\frac 32\log 2$.
Hence we have
$A_{2,1}+2A_{2,2}\sim \frac{\pi^2}{12}n^2$.

\end{itemize}
\item \underline{Control of $A_3$}. 

Notice that
$\sum_{1\le k_1< k_2<\ell_2<\ell_1\le n}\PP(E_{k_1,\ell_1}\cap E_{k_2,\ell_2})$
and
$\sum_{1\le k_1< k_2<\ell_2<\ell_1\le n}\PP(E_{k_1,\ell_1})\PP(E_{k_2,\ell_2})$ are
in $n^2\log n$. But we will see that their difference is in $n^2$.
According to \cite{FP09b} and to (\ref{TLL}), we have:
 \begin{equation}\label{A3}
A_3=c_1\!\!\!\!\!\!\!\!\!\!\sum_{1\le k_1\le k_2<\ell_2<\ell_1\le n}
    \left[\left(\sum_x \frac{e^{-\frac{\langle (\Sigma^2)^{-1}x,x\rangle}{2}\left(
         \frac 1{k_2-k_1}+\frac 1{\ell_1-\ell_2}\right)}}
   {2\pi\sqrt{\det \Sigma^2}(k_2-k_1)(\ell_1-\ell_2)}\right)
  -\frac{1}{(\ell_1-k_1)}\right]\PP(E_{k_2,\ell_2})+o(n^2)
\end{equation}
with the same notations as for (\ref{A2}).
Using again (\ref{clef}) and (\ref{TLL}), we obtain
\begin{eqnarray*}
A_3&=&o(n^{2})+c_1^2 \sum_{1\le k_1<k_2<\ell_2\le \ell_1\le n}
  \frac 1{\ell_2-k_2}\left[\frac 1{(\ell_1-k_1)-(\ell_2-k_2)}-\frac 1{(\ell_1-k_1)}\right]\\
&=&o(n^{2})+c_1^2 \sum_{1\le k_1<k_2<\ell_2\le \ell_1\le n}
  \frac {1} {(\ell_1-k_1)[(\ell_1-k_1)-(\ell_2-k_2)]}\\
&\sim&c_1^2n^2\int_{[0,1]^4}\frac{\mathbf{1}_{\{t+u+v+w<1\}}\, dt\, du\, dv\, dw}{(u+w)(u+v+w)}
=c_1^2n^2\int_{0\le u\le r\le s\le 1}\frac{(1-s)\, du\, dr\, ds}{rs}=\frac{c_1^2}2n^2.
\end{eqnarray*}
\end{itemize}%

\end{document}